\documentclass[10.5pt]{article}
\usepackage{amsmath}
\usepackage{amssymb}
\usepackage{amscd}
\usepackage{xspace}
\usepackage{verbatim}
\newcommand{\mysection}[1]{
\section{#1}\setcounter{equation}{0}}
\title{\bf Solutions of some nonlinear parabolic equations with initial blow-up}
\author{{\bf Waad Al Sayed}\quad
 {\bf Laurent V\'eron}\\[2mm]
{\small Laboratoire de Math\'ematiques et Physique Th\'eorique, }\\
{\small  Universit\'e Fran\c{c}ois Rabelais,  Tours,  FRANCE}}

\date{}
\begin{document}
\maketitle
\noindent{\small {\bf Abstract} We study the existence and uniqueness of  solutions of $\partial_tu-\Delta u+u^q=0$ ($q>1$) in $\Omega\times (0,\infty)$ where $\Omega\subset\mathbb R^N$ is a domain with a compact boundary, subject to the conditions $u=f\geq 0$ on $\partial\Omega\times (0,\infty)$ and the initial condition $\lim_{t\to 0}u(x,t)=\infty$. By means of Brezis' theory of maximal monotone operators in Hilbert spaces, we construct a minimal solution when $f=0$, whatever is the regularity of the boundary of the domain. When $\partial\Omega$ satisfies the parabolic Wiener criterion and $f$ is continuous, we construct a maximal solution and prove that it is the unique solution which blows-up at $t=0$.
}

\noindent
{\it \footnotesize 1991 Mathematics Subject Classification}. {\scriptsize
35K60}.\\
{\it \footnotesize Key words}. {\scriptsize Parabolic equations, singular solutions, semi-groups of contractions, maximal monotone operators, Wiener criterion.}
\vspace{1mm}
\hspace{.05in}

\newcommand{\txt}[1]{\;\text{ #1 }\;}
\newcommand{\tbf}{\textbf}
\newcommand{\tit}{\textit}
\newcommand{\tsc}{\textsc}
\newcommand{\trm}{\textrm}
\newcommand{\mbf}{\mathbf}
\newcommand{\mrm}{\mathrm}
\newcommand{\bsym}{\boldsymbol}
\newcommand{\scs}{\scriptstyle}
\newcommand{\sss}{\scriptscriptstyle}
\newcommand{\txts}{\textstyle}
\newcommand{\dsps}{\displaystyle}
\newcommand{\fnz}{\footnotesize}
\newcommand{\scz}{\scriptsize}
\newcommand{\be}{
\begin{equation}
}
\newcommand{\bel}[1]{
\begin{equation}
\label{#1}}
\newcommand{\ee}{
\end{equation}
}
\newcommand{\eqnl}[2]{
\begin{equation}
\label{#1}{#2}
\end{equation}
}
\newtheorem{subn}{\name}
\renewcommand{\thesubn}{}
\newcommand{\bsn}[1]{\def\name{#1}
\begin{subn}}
\newcommand{\esn}{
\end{subn}}
\newtheorem{sub}{\name}[section]
\newcommand{\dn}[1]{\def\name{#1}}   
\newcommand{\bs}{
\begin{sub}}
\newcommand{\es}{
\end{sub}}
\newcommand{\bsl}[1]{
\begin{sub}\label{#1}}
\newcommand{\bth}[1]{\def\name{Theorem}
\begin{sub}\label{t:#1}}
\newcommand{\blemma}[1]{\def\name{Lemma}
\begin{sub}\label{l:#1}}
\newcommand{\bcor}[1]{\def\name{Corollary}
\begin{sub}\label{c:#1}}
\newcommand{\bdef}[1]{\def\name{Definition}
\begin{sub}\label{d:#1}}
\newcommand{\bprop}[1]{\def\name{Proposition}
\begin{sub}\label{p:#1}}
\newcommand{\R}{\eqref}
\newcommand{\rth}[1]{Theorem~\ref{t:#1}}
\newcommand{\rlemma}[1]{Lemma~\ref{l:#1}}
\newcommand{\rcor}[1]{Corollary~\ref{c:#1}}
\newcommand{\rdef}[1]{Definition~\ref{d:#1}}
\newcommand{\rprop}[1]{Proposition~\ref{p:#1}}
\newcommand{\BA}{
\begin{array}}
\newcommand{\EA}{
\end{array}}
\newcommand{\BAN}{\renewcommand{\arraystretch}{1.2}
\setlength{\arraycolsep}{2pt}
\begin{array}}
\newcommand{\BAV}[2]{\renewcommand{\arraystretch}{#1}
\setlength{\arraycolsep}{#2}
\begin{array}}
\newcommand{\BSA}{
\begin{subarray}}
\newcommand{\ESA}{
\end{subarray}}
\newcommand{\BAL}{
\begin{aligned}}
\newcommand{\EAL}{
\end{aligned}}
\newcommand{\BALG}{
\begin{alignat}}
\newcommand{\EALG}{
\end{alignat}}
\newcommand{\BALGN}{
\begin{alignat*}}
\newcommand{\EALGN}{
\end{alignat*}}
\newcommand{\note}[1]{\textit{#1.}\hspace{2mm}}
\newcommand{\Proof}{\note{Proof}}
\newcommand{\qeda}{\hspace{10mm}\hfill $\square$}
\newcommand{\qed}{\\
${}$ \hfill $\square$}
\newcommand{\Remark}{\note{Remark}}
\newcommand{\modin}{$\,$\\
[-4mm] \indent}
\newcommand{\forevery}{\quad \forall}
\newcommand{\set}[1]{\{#1\}}
\newcommand{\setdef}[2]{\{\,#1:\,#2\,\}}
\newcommand{\setm}[2]{\{\,#1\mid #2\,\}}
\newcommand{\lra}{\longrightarrow}
\newcommand{\lla}{\longleftarrow}
\newcommand{\llra}{\longleftrightarrow}
\newcommand{\Lra}{\Longrightarrow}
\newcommand{\Lla}{\Longleftarrow}
\newcommand{\Llra}{\Longleftrightarrow}
\newcommand{\warrow}{\rightharpoonup}
\newcommand{
\paran}[1]{\left (#1 \right )}
\newcommand{\sqbr}[1]{\left [#1 \right ]}
\newcommand{\curlybr}[1]{\left \{#1 \right \}}
\newcommand{\abs}[1]{\left |#1\right |}
\newcommand{\norm}[1]{\left \|#1\right \|}
\newcommand{
\paranb}[1]{\big (#1 \big )}
\newcommand{\lsqbrb}[1]{\big [#1 \big ]}
\newcommand{\lcurlybrb}[1]{\big \{#1 \big \}}
\newcommand{\absb}[1]{\big |#1\big |}
\newcommand{\normb}[1]{\big \|#1\big \|}
\newcommand{
\paranB}[1]{\Big (#1 \Big )}
\newcommand{\absB}[1]{\Big |#1\Big |}
\newcommand{\normB}[1]{\Big \|#1\Big \|}

\newcommand{\thkl}{\rule[-.5mm]{.3mm}{3mm}}
\newcommand{\thknorm}[1]{\thkl #1 \thkl\,}
\newcommand{\trinorm}[1]{|\!|\!| #1 |\!|\!|\,}
\newcommand{\bang}[1]{\langle #1 \rangle}
\def\angb<#1>{\langle #1 \rangle}
\newcommand{\vstrut}[1]{\rule{0mm}{#1}}
\newcommand{\rec}[1]{\frac{1}{#1}}
\newcommand{\opname}[1]{\mbox{\rm #1}\,}
\newcommand{\supp}{\opname{supp}}
\newcommand{\dist}{\opname{dist}}
\newcommand{\myfrac}[2]{{\displaystyle \frac{#1}{#2} }}
\newcommand{\myint}[2]{{\displaystyle \int_{#1}^{#2}}}
\newcommand{\mysum}[2]{{\displaystyle \sum_{#1}^{#2}}}
\newcommand {\dint}{{\displaystyle \int\!\!\int}}
\newcommand{\q}{\quad}
\newcommand{\qq}{\qquad}
\newcommand{\hsp}[1]{\hspace{#1mm}}
\newcommand{\vsp}[1]{\vspace{#1mm}}
\newcommand{\ity}{\infty}
\newcommand{\prt}{\partial}
\newcommand{\sms}{\setminus}
\newcommand{\ems}{\emptyset}
\newcommand{\ti}{\times}
\newcommand{\pr}{^\prime}
\newcommand{\ppr}{^{\prime\prime}}
\newcommand{\tl}{\tilde}
\newcommand{\sbs}{\subset}
\newcommand{\sbeq}{\subseteq}
\newcommand{\nind}{\noindent}
\newcommand{\ind}{\indent}
\newcommand{\ovl}{\overline}
\newcommand{\unl}{\underline}
\newcommand{\nin}{\not\in}
\newcommand{\pfrac}[2]{\genfrac{(}{)}{}{}{#1}{#2}}

\def\ga{\alpha}     \def\gb{\beta}       \def\gg{\gamma}
\def\gc{\chi}       \def\gd{\delta}      \def\ge{\epsilon}
\def\gth{\theta}                         \def\vge{\varepsilon}
\def\gf{\phi}       \def\vgf{\varphi}    \def\gh{\eta}
\def\gi{\iota}      \def\gk{\kappa}      \def\gl{\lambda}
\def\gm{\mu}        \def\gn{\nu}         \def\gp{\pi}
\def\vgp{\varpi}    \def\gr{\rho}        \def\vgr{\varrho}
\def\gs{\sigma}     \def\vgs{\varsigma}  \def\gt{\tau}
\def\gu{\upsilon}   \def\gv{\vartheta}   \def\gw{\omega}
\def\gx{\xi}        \def\gy{\psi}        \def\gz{\zeta}
\def\Gg{\Gamma}     \def\Gd{\Delta}      \def\Gf{\Phi}
\def\Gth{\Theta}
\def\Gl{\Lambda}    \def\Gs{\Sigma}      \def\Gp{\Pi}
\def\Gw{\Omega}     \def\Gx{\Xi}         \def\Gy{\Psi}

\def\CS{{\mathcal S}}   \def\CM{{\mathcal M}}   \def\CN{{\mathcal N}}
\def\CR{{\mathcal R}}   \def\CO{{\mathcal O}}   \def\CP{{\mathcal P}}
\def\CA{{\mathcal A}}   \def\CB{{\mathcal B}}   \def\CC{{\mathcal C}}
\def\CD{{\mathcal D}}   \def\CE{{\mathcal E}}   \def\CF{{\mathcal F}}
\def\CG{{\mathcal G}}   \def\CH{{\mathcal H}}   \def\CI{{\mathcal I}}
\def\CJ{{\mathcal J}}   \def\CK{{\mathcal K}}   \def\CL{{\mathcal L}}
\def\CT{{\mathcal T}}   \def\CU{{\mathcal U}}   \def\CV{{\mathcal V}}
\def\CZ{{\mathcal Z}}   \def\CX{{\mathcal X}}   \def\CY{{\mathcal Y}}
\def\CW{{\mathcal W}} \def\CQ{{\mathcal Q}}
\def\BBA {\mathbb A}   \def\BBb {\mathbb B}    \def\BBC {\mathbb C}
\def\BBD {\mathbb D}   \def\BBE {\mathbb E}    \def\BBF {\mathbb F}
\def\BBG {\mathbb G}   \def\BBH {\mathbb H}    \def\BBI {\mathbb I}
\def\BBJ {\mathbb J}   \def\BBK {\mathbb K}    \def\BBL {\mathbb L}
\def\BBM {\mathbb M}   \def\BBN {\mathbb N}    \def\BBO {\mathbb O}
\def\BBP {\mathbb P}   \def\BBR {\mathbb R}    \def\BBS {\mathbb S}
\def\BBT {\mathbb T}   \def\BBU {\mathbb U}    \def\BBV {\mathbb V}
\def\BBW {\mathbb W}   \def\BBX {\mathbb X}    \def\BBY {\mathbb Y}
\def\BBZ {\mathbb Z}

\def\GTA {\mathfrak A}   \def\GTB {\mathfrak B}    \def\GTC {\mathfrak C}
\def\GTD {\mathfrak D}   \def\GTE {\mathfrak E}    \def\GTF {\mathfrak F}
\def\GTG {\mathfrak G}   \def\GTH {\mathfrak H}    \def\GTI {\mathfrak I}
\def\GTJ {\mathfrak J}   \def\GTK {\mathfrak K}    \def\GTL {\mathfrak L}
\def\GTM {\mathfrak M}   \def\GTN {\mathfrak N}    \def\GTO {\mathfrak O}
\def\GTP {\mathfrak P}   \def\GTR {\mathfrak R}    \def\GTS {\mathfrak S}
\def\GTT {\mathfrak T}   \def\GTU {\mathfrak U}    \def\GTV {\mathfrak V}
\def\GTW {\mathfrak W}   \def\GTX {\mathfrak X}    \def\GTY {\mathfrak Y}
\def\GTZ {\mathfrak Z}   \def\GTQ {\mathfrak Q}

\font\Sym= msam10 
\def\SYM#1{\hbox{\Sym #1}}
\newcommand{\bdw}{\prt\Gw\xspace}
\medskip
\mysection {Introduction}

\setcounter{equation}{0}
Let $\Gw$ be a domain of $\BBR^N$ ($N\geq 1$) with a compact boundary, $Q_{\infty}^\Gw=\Gw\ti (0,\infty)$ and $q>1$. This article deals with the question of the solvability of the following Cauchy-Dirichlet  problem $\CP^{\Gw,f}$
\begin{equation}\label{Q1}\left\{\BA {l}
\prt_{t}u-\Gd u+|u|^{q-1}u=0\quad\text{in }\,Q_{\infty}^\Gw\\[2mm]
\phantom{\prt_{t}u-\Gd u+|u|^{q-1}}
u=f\quad\text{on }\,\prt \Gw\ti (0,\infty)\\[2mm]
\phantom{\prt_{t}u-}
\lim_{t\to 0}u(x,t)=\infty\qquad\forall x\in\Gw.
  \EA\right.\end{equation}
 If no assumption of regularity is made on $\prt\Gw$, the boundary data $u=f$ cannot be prescribed in sense of continuous functions. However, the case $f=0$ can be treated if the vanishing condition on $\prt \Gw\ti (0,\infty)$ is understood in the $H^1_{0}$ local sense. We construct a positive solution $\underline u_{\Gw}$ of (\ref{Q1}) with $f=0$ belonging to $C(0,\infty; H^1_{0}(\Gw)\cap L^{q+1}(\Gw))$ thanks to Brezis results of contractions semigroups generated by subdifferential of proper convex functions in Hilbert spaces. We can also consider an internal increasing approximation of $\Gw$ by smooth bounded domains $\Gw^n$ such that $\Gw=\cup_{n}\Gw^n$. For each of these domains, there exists a maximal solution $\overline u_{\Gw^n}$ of problem $\CP^{\Gw^n,0}$.
Furthermore the sequence $\{\overline u_{\Gw^n}\}$ is increasing. The limit function $ u_{\Gw}:=\lim_{n\to\infty}\overline u_{\Gw^n}$ is the natural candidate to be the minimal positive solution of  a solution of $\CP^{\Gw,0}$. We prove that $\underline u_{\Gw}= u_{\Gw}$. If $\prt\Gw$ satisfies the parabolic Wiener criterion \cite {Zi}, there truly exist solutions of $\CP^{\Gw,0}$. We construct a maximal solution $\overline u_{\Gw}$ of this problem. Our main result is the following:\medskip
  
\noindent{\bf Theorem 1. }{\it If $\prt\Gw$ is compact and satisfies the parabolic Wiener criterion, there holds}
 $$\overline u_{\Gw}= \underline u_{\Gw}.$$
 
 In the last section, we consider the full problem $\CP^{\Gw,f}$. Under the same regularity and boundedness assumption on $\prt\Gw$ we construct a maximal solution $\overline u_{\Gw,f}$ and we prove
 \medskip
  
\noindent{\bf Theorem 2. }{\it If $\prt\Gw$ is compact and satisfies the parabolic Wiener criterion, and if $f\in C(0,\infty;\prt\Gw)$ is nonnegative, $\overline u_{\Gw,f}$ is the only positive solution to problem }$\CP^{\Gw,f}$.\medskip

These type of results are to be compared with the ones obtained by the same authors \cite{AV} in which paper the following problem is considered
 \begin{equation}\label{Q3}\left\{\BA {l}
\,\,\;\;\prt_{t}u-\Gd u+|u|^{q-1}u=0\quad\text{in }\,Q_{\infty}^\Gw\\[2mm]

\lim_{\dist (x,\prt\Gw)\to 0}u(x,t)=\infty\quad\text{locally uniformly on } (0,\infty)\\[2mm]
\phantom{\prt_{t}u,-----}
u(x,0)=f\qquad\forall x\in\Gw.
  \EA\right.\end{equation}
  
  In the above mentioned paper, it is proved two types of existence and uniqueness result with $f\in L^1_{loc}(\Gw)$, $f\geq 0$: either if $\prt\Gw=\prt\overline\Gw^c$ and $1<q<N/(N-2)$, or if $\prt\Gw$ is locally the graph of a continuous function and $q>1$.\medskip
  
  Our paper is organized as follows: 1- Introduction. 2- Minimal and maximal solutions. 3- Uniqueness of large solutions. 4- Bibliography.
 \mysection {Minimal and maximal solutions}

\setcounter{equation}{0}
Let $q>1$ and $\Gw$ be a proper domain of $\BBR^N$, $N>1$ with a non-empty compact boundary. We set $Q_{\infty}^{\Gw}=\Gw\ti(0,\infty)$ and consider the following problem
\begin{equation}\label{MS0}\left\{\BA{l}Ê 
\prt_{t}u-\Gd u+u^q=0\quad\text{in } \Gw\ti (0,\infty)\\
\phantom{-,,--}
\!u(x,t)=0\quad\quad\text{on } \prt\Gw\ti (0,\infty).
\EA\right.\end{equation}
If there is no regularity assumption on $\prt\Gw$, a natural way to consider the boundary condition is to impose $u(.,t)\in H^1_{0}(\Gw)$. The Hilbertian framework for this equation has been studied by Brezis in a key article \cite{Br1} (see also the monography \cite {Br2} for a full treatment of related questions) in considering the maximal monotone operator $v\mapsto A(v):=-\Gd v+|v|^{q-1}v$ seen as the subdifferential of the proper lower semi-continuous function
\begin{equation}\label{MS1}J_{\Gw}(v)=\left\{\BA{c}Ê 
\myint{\Gw}{}\left(\myfrac{1}{2}|\nabla v|^2+\myfrac{1}{q+1}|v|^{q+1}\right) dx\quad\text{if }v\in H^1_{0}(\Gw)\cap L^{q+1}(\Gw)\\[4mm]\phantom{--------}
\infty\phantom{-------,}\text{if }v\notin H^1_{0}(\Gw)\cap L^{q+1}(\Gw).
\EA\right.\end{equation}
In that case, the domain of $A=\prt J_{\Gw}$ is $D(A):=\{u\in H^1_{0}(\Gw)\cap L^{q+1}(\Gw):\Gd u\in L^2(\Gw)\}$, and we endow $D_{\Gw}(-\Gd,)$ with the graph norm of the Laplacian in $ H^1_{0}(\Gw)$
$$\norm v_{D_{\Gw}(-\Gd)}=\left(\myint{\Gw}{}\left((\Gd v)^2+|\nabla v|^2+v^2\right)dx\right)^{1/2}.
$$
Brezis' result is the following.

\bth{Br}Given $u_{0}\in L^2(\Gw)$ there exists a unique function $v\in L_{loc}^2(0,\infty;D_{\Gw}(-\Gd))\cap C(0,\infty; H^1_{0}(\Gw)\cap L^{q+1}(\Gw))$ such that $\prt_tv\in L_{loc}^2(0,\infty;L^2(\Gw))$ satisfying
\begin{equation}\label{HB}\left\{\BA{l}Ê 
\prt_{t}v-\Gd v+|v|^{q-1}v=0\quad \text{a.e. in }Q_{\infty}^{\Gw}\\
\phantom{---.--}
v(.,0)=u_{0}\quad \text{a.e. in }\Gw.
\EA\right.\end{equation}
Furthermore the mapping $(t,u_{0})\mapsto v(t,.)$ defines an order preserving contraction semigroup in $L^2(\Gw)$, denoted by $S^{\prt J_{\Gw}}(t)[u_{0}]$, and the following estimate holds
\begin{equation}\label{HB'}
\norm{\prt_{t}v(t,.)}_{L^2(\Gw)}\leq\myfrac{1}{t\sqrt 2}\norm{u_{0}}_{L^2(\Gw)}.
\end{equation}
\es

From this result, we have only to consider solutions of (\ref{MS0}) with the above regularity. 

\bdef{I(Q)}We denote by $\CI(Q_{\infty}^{\Gw})$ the set of positive functions $u\in L_{loc}^2(0,\infty;D_{\Gw}(-\Gd))\cap C(0,\infty; H^1_{0}(\Gw)\cap L^{q+1}(\Gw))$ such that $\prt_{t}u\in L_{loc}^2(0,\infty;L^2(\Gw))$ satisfying 
\begin{equation}\label{E1}
\prt_{t}u-\Gd u+|u|^{q-1}u=0
\end{equation} 
in the semigroup sense, i. e.
\begin{equation}\label{HB1}\myfrac{du}{dt}+\prt J_{\Gw}(u)=0\quad\text {a.e. in }(0,\infty).
\end{equation}
\es
If $\Gw$ is not bounded it is usefull to introduce another class which takes into account the Dirichlet condition on $\prt\Gw$: we assume that $\Gw^c\subset B_{R_{0}}$, denote by $\Gw_{R}=\Gw\cap B_{R}$ ($R\geq R_{0}$) and by
$\tilde H_{0}^1(\Gw_{R})$ the closure in $H^1_{0}(\Gw_{R})$ of the restrictions to $\Gw_{R}$ of functions in $C^\infty_{0}(\Gw)$, thus we endow $D_{\Gw_{R}}(-\Gd,)$ with the graph  norm of the Laplacian in $ \tilde H^1_{0}(\Gw_{R})$
$$\norm v_{D_{\Gw_{R}}(-\Gd)}=\left(\myint{\Gw_{R}}{}\left((\Gd v)^2+|\nabla v|^2+v^2\right)dx\right)^{1/2}.
$$

\bdef{I(Qloc)} If $\Gw$ is not bounded but $\Gw^c\subset B_{R_{0}}$, we denote by  $\CI(Q_{\infty}^{\Gw_{loc}})$ the set of positive functions $u\in L^2_{loc}(Q^\Gw_{\infty})$ such that, for any $R>R_{0}$,
$u\in L_{loc}^2(0,\infty;D_{\Gw_{R}}(-\Gd))\cap C(0,\infty; \tilde H^1_{0}(\Gw_{R})\cap L^{q+1}(\Gw_{R}))$, $\prt_{t}u\in L_{loc}^2(0,\infty;L^2(\Gw_{R}))$ and $u$ satisfies (\ref{E1})  in a. e. in $Q^\Gw_{\infty}$.
\es

 \blemma{ext}If $u\in \CI(Q_{\infty}^{\Gw})$ or $\CI(Q_{\infty}^{\Gw_{loc}})$, its extension $\tilde u$ by zero outside $\Gw$ is a subsolution of (\ref{MS0}) in $(0,\infty)\ti\BBR^N$ such that $\tilde u\in C(0,\infty;H^1_{0}(\BBR^N)\cap L^{q+1}(\BBR^N))$ and $\prt_t\tilde u\in L^2_{loc}(0,\infty;L^2(\BBR^N))$.
  \es
  \Proof The proof being similar in the two cases, we assume $\Gw$ bounded. We first notice that $\tilde u\in C(0,\infty;H^1_{0}(\BBR^N))$ since $\norm {\tilde u}_{H^1_{0}(\BBR^N)}=\norm { u}_{H^1_{0}(\Gw)}$. For $\gd>0$ we set
  $$
  P_{\gd}(r)=\left\{ \BA {l}r-3\gd/2\qquad\qquad\text {if }r\geq 2\gd\\
  r^2/2\gd-r+\gd/2\quad\text {if }\gd<r< 2\gd\\
  0\qquad\quad\qquad\qquad\text {if }r\leq \gd  \EA
  \right.
$$
and denote by $u_{\gd}$ the extension of $P_{\gd}(u)$ by zero outside $Q_{\infty}^{\Gw}$. Since $ u_{\gd\,t}=P'_{\gd}(u) \prt_tu$, then
$ u_{\gd\,t}\in L^2_{loc} (0,\infty; L^2(\BBR^N))$ and $\norm{u_{\gd\,t}}_{L^2}\leq \norm{\prt_tu}_{L^2}$. In the same way $\nabla u_{\gd}=P'_{\gd}(u)\nabla u$, thus $ u_{\gd}\in L^2_{loc} (0,\infty; H^1_{0}(\BBR^N))$ and $\norm{u_{\gd}}_{H^1_{0}}\leq \norm{u}_{H^1_{0}}$. Finally $-\Gd u_{\gd}=-P'_{\gd}(u)\Gd u-P''_{\gd}(u)\abs{\nabla u}^2.$ Using the fact that $P'_{\gd}u^q\geq u_{\gd}^q$, we derive from
(\ref{HB1})
$$\prt_{t}u_{\gd}-\Gd u_{\gd}+u^q_{\gd}\leq 0
$$
in the sense that
\begin{equation}\label{HB2}
\dint_{Q^{\BBR^N}_\infty}\left(\prt_{t}u_{\gd}\gz+\nabla u_{\gd}.\nabla\gz+u^q_{\gd}\gz\right)dxdt\leq 0
\end{equation}
for all $\gz\in C^\infty((0,\infty)\ti\BBR^N)$, $\gz\geq 0$. Actually, $C^\infty((0,\infty)\ti\BBR^N)$ can be replaced by $L^2(\ge,\infty;H^1_{0}(\BBR^N))\cap L^{q'}((\ge,\infty)\ti\BBR^N)$. Letting $\gd\to 0$ and using Fatou's theorem implies that (\ref{HB2}) holds with $u_{\gd}$ replaced by $\tilde u$.
\qeda
 \blemma{sup} For any $u\in \CI(Q_{\infty}^{\Gw})$, there holds
\begin{equation}\label{MS3}
u(x,t)\leq \left(\myfrac{1}{(q-1)t}\right)^{1/(q-1)}:=\phi_{q}(t)\quad\forall (x,t)\in Q_{\infty}^{\Gw}.
\end{equation}
\es
\Proof Let $\gt>0$. Since the function $\phi_{q,\gt}$ defined by $\phi_{q,\gt}(t)=\phi_{q}(t-\gt)$ is a solution of 
$$\phi'_{q,\gt}+\phi^q_{q,\gt}=0
$$
and $(u-\phi_{q,\gt})_{+}\in C(0,\infty;H^1_{0}(\Gw))$, there holds
$$\myfrac{1}{2}\myfrac{d}{dt}\myint{\Gw}{}(u-\phi_{q,\gt})^2_{+}dx
+\dint_{\!\!Q_{\infty}^{\Gw}}\left(\nabla u.\nabla (u-\phi_{q,\gt})_{+}+
(u^q-\phi^q_{q,\gt})(u-\phi_{q,\gt})_{+}\right)dxdt=0.
$$
Thus $s\mapsto \norm{ (u-\phi_{q,\gt})_{+}(s)}_{L^2}$
is nonincreasing. By Lebesgue's theorem, 
$$\lim_{s\downarrow\gt}\norm{ (u-\phi_{q,\gt})_{+}(s)}_{L^2}=0,$$
thus $u(x,t)\leq \phi_{q,\gt}(t)$ a.e. in $\Gw$. Letting $\gt\downarrow 0$ and using the continuity yields to (\ref{MS3}).\qeda

\bth{max} For any $q>1$, the set $\CI(Q_{\infty}^{\Gw})$ admits a least upper bound $\underline u_{\Gw}$ for the order relation. If $\Gw$ is bounded, $\underline u_{\Gw}\in \CI(Q_{\infty}^{\Gw})$; if it is not the case, then $\underline u_{\Gw}\in \CI(Q_{\infty}^{\Gw_{loc}})$.
\es
\Proof \noindent {\it Step 1- Construction of $\underline u_{\Gw}$ when $\Gw$ is bounded.} For $k\in\BBN^*$ we consider the solution $v=v_{k}$ (in the sense of \rth{Br} with the corresponding maximal operator in $L^2(\Gw)$) of
\begin{equation}\label{I1}\left\{\BA{l}Ê 
\prt_{t}v-\Gd v+v^q=0\qquad\text{in } \Gw\ti (0,\infty)\\
\phantom{-,,--}
\!v\!(x,t)=0\quad\quad\text{in } \prt\Gw\ti (0,\infty)\\
\phantom{-,--}
v(x,0)=k\quad\quad\text{in } \Gw.
\EA\right.\end{equation}
When $k\to\infty$, $v_{k}$ increases and converges to some $\underline u_{\Gw}$. Because of (\ref{MS3}) and the fact that $\Gw$ is bounded, $\underline u_{\Gw}(t,.)\in L^2(\Gw)$ for $t>0$. It follows from the closedness of maximal monotone operators that $\underline u_{\Gw}\in L^2_{loc}(0,\infty;D_{\Gw}(-\Gd))\cap C(0,\infty; H^1_{0}(\Gw)\cap L^{q+1}(\Gw))$, $\prt_{t}\underline u_{\Gw}\in L_{loc}^2(0,\infty;L^2(\Gw))$ and 
\begin{equation}\label{I2}\myfrac{d\underline u_{\Gw}}{dt}+\prt J_{\Gw}(\underline u_{\Gw})=0\quad\text {a.e. in }(0,\infty).
\end{equation}
Thus $\underline u_{\Gw}\in\CI(Q_{\infty}^\Gw)$. For $\gt,\ge>0$, the function $t\mapsto\underline u_{\Gw}(x,t-\gt)+\ge$ is a supersolution of (\ref{MS0}). Let $u\in \CI(Q^\Omega_\infty)$; for $k>\phi_{q}(\gt)$, the function $(x,t)\mapsto(u(x,t)-\underline u_{\Gw}(x;t-\gt)-\ge)_{+}$ is a subsolution of (\ref{MS0}) and belongs to $C(\gt,\infty;H^1_{0}(\Gw))$. Since it vanishes at $t=\gt$, it follows from Brezis' result that it is identically zero, thus $u(x,t)\leq\underline u_{\Gw}(x,t-\gt)+\ge$. Letting $\ge,\gt\downarrow 0$ implies the claim.\smallskip

\noindent {\it Step 2- Construction of $\underline u_{\Gw}$ when $\Gw$ is unbounded.} We assume that $\prt\Gw\subset B_{R_{0}}$ and for $n>R_{0}$, we recall that $\Gw_{n}=\Gw\cap B_{n}$. For $k>0$, we denote by $\underline u_{\Gw_{n}}$ the solution obtained in Step 1. Then $\underline u_{\Gw_{n}}=\lim_{k\to\infty}v_{n,k}$ where $v_{n,k}$ is the solution, in the sense of maximal operators in $\Gw_{n}$ of
\begin{equation}\label{I3}\left\{\BA {l}\myfrac{dv_{n,k}}{dt}+\prt J_{\Gw_{n}}(v_{n,k})=0\quad\text {a.e. in }(0,\infty)
\\\phantom{\myfrac{dv_{n,k}}{dt}+\prt J_{\Gw_{n}}}
v_{n,k}(0)=k.
\EA\right.\end{equation}
It follows from \rlemma{ext} that the extension $\tilde v_{n,k}$ by $0$ of $v_{n,k}$ in 
$\Gw_{n+1}$ is a subsolution for the equation satisfied by $v_{n+1,k}$, with a smaller initial data, therefore $\tilde v_{n,k}\leq v_{n+1,k}$. This implies $\tilde {\underline u}_{\Gw_{n}}\leq u_{\Gw_{n+1}}$. Thus we define $\underline u_{\Gw}=\lim_{n\to\infty}\underline u_{\Gw_{n}}$. It follows from \rlemma {sup} and standard regularity results for parabolic equations  that $u=\underline u_{\Gw}$ satisfies
\begin{equation}\label{Equ}
\prt_{t}u-\Gd u+u^q=0
\end{equation}
 in $Q_{\infty}^{\Gw}$. Multiplying 
\begin{equation}\label{Equ-n}\myfrac{d\underline u_{\Gw_{n}}}{dt}+\prt J_{\Gw_{n}}(\underline u_{\Gw_{n}})=0
\end{equation}
by $\eta^2 \underline u_{\Gw_{n}}$ where $\eta\in C^\infty_0(\BBR^N)$ and integrating over $\Gw_{n}$, yields to
$$2^{-1}\myfrac{d}{dt}\myint{\Gw_{n}}{}\eta^2 \underline u^2_{\Gw_{n}}dx+\myint{\Gw_{n}}{}\left(|\nabla \underline u_{\Gw_{n}}|^2+
\underline u^{q+1}_{\Gw_{n}}\right)\eta^2dx+
2\myint{\Gw_{n}}{}\nabla \underline u_{\Gw_{n}}.\nabla\eta\,\eta \underline u_{\Gw_{n}}dx=0.
$$
Thus, by Young's inequality,
$$2^{-1}\myfrac{d}{dt}\myint{\Gw_{n}}{}\eta^2\underline u^2_{\Gw_{n}}dx+\myint{\Gw_{n}}{}\left(2^{-1}|\nabla \underline u_{\Gw_{n}}|^2+\underline
u^{q+1}_{\Gw_{n}}\right)\eta^2dx\leq 
2\myint{\Gw_{n}}{}|\nabla\eta|^2 \underline u^2_{\Gw_{n}}dx.
$$
If we assume that $0\leq \eta\leq 1$, $\eta=1$ on $B_{R}$ ($R>R_{0}$) and $\eta=0$ on $B_{2R}^c$, we derive, for any $0<\gt<t$,
\begin{equation}\label{I5}\BA {l}
2^{-1}\myint{\Gw_{n}}{} \underline u_{\Gw_{n}}^2(.,t)\eta^2dx+\myint{\gt}{t}\myint{\Gw_{n}}{}\left(2^{-1}\abs{\nabla  \underline u_{\Gw_{n}}}^2 +\underline u_{\Gw_{n}}^{q+1} \right)\eta^2dxds\\
\phantom{--------------}
\leq 2\myint{\gt}{t}\myint{\Gw_{n}}{} {\underline u}_{\Gw_{n}}^2|\nabla \eta|^2dxds
+2^{-1}\myint{\Gw_{n}}{}{\underline u}_{\Gw_{n}}^2(.,\gt)\eta^2dx.
\EA\end{equation}
From this follows, if $n>2R$,
\begin{equation}\label{I6}\BA {l}
2^{-1}\myint{\Gw\cap B_{R}}{}\underline u_{\Gw_{n}}^2(.,t)dx+\myint{\gt}{t}\myint{\Gw\cap B_{R}}{}\left(2^{-1}\abs{\nabla \underline u_{\Gw_{n}}}^2 +\underline u_{\Gw_{n}}^{q+1} \right) dxds\leq CR^N(t+1)\gt^{-2/(q-1)}.
\EA\end{equation}
If we let $n\to\infty$ we derive by Fatou's lemma
\begin{equation}\label{I7}\BA {l}
2^{-1}\myint{\Gw\cap B_{R}}{}\underline u_{\Gw}^2(.,t)dx+\myint{\gt}{t}\myint{\Gw\cap B_{R}}{}\left(2^{-1}\abs{\nabla\underline u_{\Gw}}^2 +\underline u_{\Gw}^{q+1} \right) dxds\leq CR^N(t+1)\gt^{-2/(q-1)}.
\EA\end{equation}
For $\gt>0$ fixed, we multiply (\ref{Equ-n}) by $(t-\gt)\eta^2d\underline u_{\Gw_{n}}/dt $, integrate  on $(\gt,t)\ti\Gw_{n}$ and get
$$\BA {l}(t-\gt)\myint{\Gw_{n}}{}\abs{\myfrac{d\underline u_{\Gw_{n}}}{dt}}^2\eta^2dx
+\myfrac{d}{dt}(t-\gt)\myint{\Gw_{n}}{}\left(\myfrac{|\nabla{\underline u_{\Gw_{n}}}|^2}{2}+\myfrac{\underline u_{\Gw_{n}}^{q+1}}{q+1}\right)\eta^2dx
\\\phantom{------}
=\myint{\Gw_{n}}{}\left(\myfrac{|\nabla{\underline u_{\Gw_{n}}}|^2}{2}+\myfrac{\underline u_{\Gw_{n}}^{q+1}}{q+1}\right)\eta^2dx-2(t-\gt)\myint{\Gw_{n}}{}\nabla \underline u_{\Gw_{n}}.\nabla\eta\myfrac{d\underline u_{\Gw_{n}}}{dt}\eta dx.
\EA$$
Since
$$2(t-\gt)\abs{\myint{\Gw_{n}}{}\nabla \underline u_{\Gw_{n}}.\nabla\eta\myfrac{d\underline u_{\Gw_{n}}}{dt}\eta dx}
\leq \myfrac{(t-\gt)}{2}\myint{\Gw_{n}}{}\abs{\myfrac{d\underline u_{\Gw_{n}}}{dt}}^2\eta^2dx
+4(t-\gt)\myint{\Gw_{n}}{}\abs{\nabla \underline u_{\Gw_{n}}}^2\abs{\nabla \eta}^2dx,
$$
we get, in assuming again $n>2R$,
\begin{equation}\label{I8}\BA {l}
2^{-1}\myint{\gt}{t}\myint{\Gw}{}(s-\gt)\abs{\myfrac{d\underline u_{\Gw_{n}}}{dt}}^2\eta^2dxds+(t-\gt)\myint{\Gw}{}\left(\myfrac{|\nabla{\underline u_{\Gw_{n}}}|^2}{2}+\myfrac{\underline u_{\Gw_{n}}^{q+1}}{q+1}\right)\eta^2dx\\
\phantom{-------------------}
\leq 4\myint{\gt}{t}(s-\gt)\myint{\Gw}{}\abs{\nabla \underline u_{\Gw_{n}}}^2\abs{\nabla\eta }^2dxds,
\EA\end{equation}
from which follows,
\begin{equation}\label{I9}\BA {l}
2^{-1}\myint{\gt}{t}\myint{\Gw\cap B_{R}}{}(s-\gt)\abs{\myfrac{d\underline u_{\Gw_{n}}}{dt}}^2 dxds+(t-\gt)\myint{\Gw\cap B_{R}}{}\left(\myfrac{|\nabla{\underline u_{\Gw_{n}}}|^2}{2}+\myfrac{\underline u_{\Gw_{n}}^{q+1}}{q+1}\right)dx\\[4mm]
\phantom{-------------------}
\leq 4\myint{\gt}{t}(s-\gt)\myint{\Gw\cap B_{2R}}{}\abs{\nabla \underline u_{\Gw_{n}}}^2 dxds.
\EA\end{equation}
The right-hand side of (\ref{I9}) remains uniformly bounded by $8C(2R)^N(t-\gt)t\gt^{-2/(q-1)}$ from (\ref{I6}). Then
\begin{equation}\label{I10}\BA {l}
2^{-1}\myint{\gt}{t}\myint{\Gw\cap B_{R}}{}(s-\gt)\abs{\myfrac{d\underline u_{\Gw_{n}}}{dt}}^2 dxds+(t-\gt)\myint{\Gw\cap B_{R}}{}\left(\myfrac{|\nabla{\underline u_{\Gw_{n}}}|^2}{2}+\myfrac{\underline u_{\Gw_{n}}^{q+1}}{q+1}\right)dx\\[4mm]
\phantom{----------------------}
\leq 8C(2R)^N(t-\gt)t\gt^{-2/(q-1)}
\EA\end{equation}
By Fatou's lemma the same estimate holds if $\underline u_{\Gw_{n}}$ is replaced by 
$\underline u_{\Gw}$. Notice also that this estimate implies that $\underline u_{\Gw}$ vanishes in the $H^1_{0}$-sense on $\prt\Gw$ since $\eta\underline u_{\Gw}\in H^1_{0}(\Gw)$ where the function $\eta\in C^\infty_{0}(\BBR^N)$ has value $1$ in $B_{R}$ and $\Gw^c\subset B_{R}$. Moreover  estimates (\ref{I7}) and (\ref{I10}) imply that $\underline u_{\Gw}$ satisfies (\ref{Equ}) a.e.,  and thus it belongs to $\CI(Q_{\infty}^{\Gw_{loc}})$.\smallskip
 
 \noindent {\it Step 3- Comparison.} At end, assume $u\in \CI(Q_{\infty}^{\Gw})$. For $R>n_{0}$ let $W_{R}$ be the maximal solution of 
\begin{equation}\label{maxell}
-\Gd W_{R}+W_{R}^q=0\quad\text{in }B_{R}.
\end{equation}
Existence follows from Keller-Osserman's construction \cite{Ke},\cite{Oss}, and the following scaling and blow-up estimates holds
\begin{equation}\label{maxell1}
W_{R}(x)=R^{-2/(q-1)}W_{1}(x/R),
\end{equation}
and
\begin{equation}\label{maxell2}
W_{R}(x)=C_{q}(R-\abs x)^{-2/(q-1)}(1+\circ (1))\,\,\text {as } |x|\to R.
\end{equation}
For $\gt>0$ set $v(x,t)=u(x,t)-\underline u_{\Gw}(x,t-\gt)-W_{R}(x)$. Then $v_{+}$ is a subsolution. Since $v(.,\gt)\in L^2(\Gw)$, $\lim_{s\downarrow\gt}||v_{+}(.,s)||_{L^2}=0$. Because $\eta\underline u_{\Gw}\in H^1_{0}(\Gw)$ for $\eta$ as above, $\eta v_{+}\in H^1_{0}(\Gw)$. Next, $supp\, v_{+}\subset \Gw\cap B_{R}$. Since $u,\underline u_{\Gw}$ are locally in $H^1$, we can always assume that their restrictions to $\prt B_{R}\ti [0,T]$ are integrable for the corresponding Hausdorff measure. Therefore Green's formula is valid, which implies
$$-\myint{\gt}{t}\myint{\Gw\cap B_{R}}{}\Gd v _+dxdt
=\myint{\gt}{t}\myint{\Gw\cap B_{R}}{}|\nabla v_{+}|^2dxdt\quad\forall t>\gt.
$$
Therefore
$$\myint{\Gw\cap B_{R}}{}v_{+}^2(x,t)dx+\myint{\gt}{t}\myint{\Gw\cap B_{R}}{}\!\!\!\!\!\left(|\nabla v_{+}|^2+(u-(\underline u_{\Gw}(.,t-\gt)+W_{R})^q)v_{+}\right)dxdt\leq \myint{\Gw\cap B_{R}}{}v_{+}^2(x,s)dx.
$$
We let $s\downarrow \gt$ and get $v_{+}=0$, equivalently $u(x,t)\leq\underline u_{\Gw}(x,t-\gt)+W_{R}(x)$. Then we let $R\to\infty$ and $\gt\to 0$ and obtain $u(x,t)\leq \underline u_{\Gw}(x,t)$, which is the claim.\qeda\medskip


\bcor{ord} Assume $\Gw^{1}\subset \Gw^{2}\subset\BBR^N$ are open domains, then $ \underline u_{\Gw^{1}}\leq\underline u_{\Gw^{2}}$. Furthermore, if  $\Gw=\cup\Gw^{n}$ where $\Gw^{n}\subset\Gw^{n+1}$, then
\begin{equation}\label{comp}
\lim_{n\to\infty}\underline u_{\Gw^{n}}=\underline u_{\Gw},
\end{equation} 
locally uniformly in $Q_{\infty}^\Gw$.\es
\Proof The first assertion follows from the proof of \rth{max}. It implies 
$$\lim_{n\to\infty}\underline u_{\Gw^{n}}=u^*_{\Gw}\leq \underline u_{\Gw},
$$
and $u^*_{\Gw}$ is a positive solution of (\ref{E1}) in $Q_{\infty}^\Gw$.
There exists a sequence $\{u_{0,m}\}\subset L^2(\Gw)$ such that 
$S^{\prt J_\Gw}(t)[u_{0,m}]\uparrow \underline u_{\Gw}$ as $n\to\infty$, locally uniformly in $Q_{\infty}^\Gw$. Set $u_{0,m,n}=u_{0,m}\chi_{_{\Gw^{n}}}$; since 
$u_{0,m,n}\to u_{0,m}$ in $L^2(\Gw)$ then $S^{\prt J_\Gw}(.)[u_{0,m,n}]\uparrow S^{\prt J_\Gw}(.)[u_{0,m}]$  in $L^\infty(0,\infty;L^2(\Gw))$. If  $\tilde v_{m,n}$ is the extension of $v_{m,n}:=S^{\prt J_{\Gw^{n}}}(.)[u_{0,m,n}]$ by zero outside $Q_{\infty}^{\Gw_{n}}$ it is a subsolution smaller than $S^{\prt J_\Gw}(.)[u_{0,m,n}]$ and $n\mapsto \tilde v_{m,n}$ is increasing; we denote by $\tilde v_{m}$ its limit as $n\to\infty$. Since for any $\gz\in C_{0}^{2,1}([0,\infty)\ti\Gw)$ we have, for $n$ large enough and $s>0$,
$$
-\myint{0}{s}\myint{\Gw}{}\left(\tilde v_{m,n}\left(\prt_{t}\gz+\Gd\gz\right)\right)dxdt=\myint{\Gw}{}u_{0,m,n}\gz(x,0)dx-\myint{\Gw}{}\tilde v_{m,n}(x,s)\gz(x,t)dx,
$$
it follows
$$
-\myint{0}{s}\myint{\Gw}{}\left(\tilde v_{m}\left(\prt_{t}\gz+\Gd\gz\right)\right)dxdt=\myint{\Gw}{}u_{0,m}\gz(x,0)dx-\myint{\Gw}{}\tilde v_{m}(x,s)\gz(x,t)dx.
$$
Clearly  $\tilde v_{m}$ is a solution of (\ref{E1}) in $Q_{\infty}^{\Gw_{n}}$. Furthermore
$$\lim_{t\to 0}\tilde v_{m}(t,.)=u_{0,m}\quad\text{a.e. in }\Gw.
$$
Because 
$$\norm{\tilde v_{m}(t,.)-u_{0,m}}_{L^2(\Gw)}\leq 2\norm{u_{0,m}}_{L^2(\Gw)},
$$
it follows from Lebesgue's theorem that $t\mapsto \tilde v_{m}(t,.)$ is continuous in $L^2(\Gw)$ at $t=0$. Furthermore, for any $t>0$ and $h\in(-t,t)$, we have from \ref{HB'},

\begin{equation}\label{ord1}\BA {l}
\norm{\tilde v_{m,n}(t+h,.)-\tilde v_{m,n}(t,.)}_{L^2(\Gw ^{n})}
\leq\myfrac{|h|}{t\sqrt 2}\norm{u_{0,m,n}}_{L^2(\Gw ^{n})}\\
\phantom{----------}
\Longrightarrow\norm{\tilde v_{m}(t+h,.)-\tilde v_{m}(t,.)}_{L^2(\Gw)}
\leq\myfrac{|h|}{t\sqrt 2}\norm{u_{0,m}}_{L^2(\Gw)}.
\EA\end{equation} 
Thus $\tilde v_{m}\in C([0,\infty);L^2(\Gw)$. By the contraction principle, $\tilde v_{m}=S^{\prt J_\Gw}(t)[u_{0,m}]$ is the unique generalized solution to (\ref{HB}). Finally, there exists an increasing sequence  $\{u_{0,m}\}\subset L^2(\Gw)$ such that for any $\ge>0$, and $\gt>0$, 
$$0<\underline u_{\Gw}-S^{\prt J_\Gw}(t)[u_{0,m}]\leq \ge/2$$
on $[\gt,\infty)\ti\Gw$. For any $m$, there exists $n_{m}$ such that
$$0<S^{\prt J_\Gw}(t)[u_{0,m}]-\tilde v_{m,n}\leq \ge/2$$
Therefore
$$0<\underline u_{\Gw}-\underline u_{\Gw^{n}}\leq \ge,$$
on $[\gt,\infty)\ti\Gw_{n}$. This implies (\ref{comp}).
\qeda

\medskip

We can also construct a minimal solution with conditional initial blow-up in the following way. Assuming that $\Gw=\cup\Gw^{m}$ where $\Gw^{m}$ are smooth bounded domains and $\overline{\Gw^{m}}\subset\Gw^{m+1}$. We denote by $u_{m}$ the solution of 
\begin{equation}\label{min1}\left\{\BA {l}
\prt_{t}u_{m}-\Gd u_{m}+|u_{m}|^{q-1}u_{m}=0\quad\text{in }Q_{\infty}^{\Gw^{m}}\\[2mm]
\phantom{\prt_{t}u_{m}-\Gd u_{m,k}+|u_{m}|^{q-1}}
u_{m}=0\quad\text{in }\prt\Gw^{m}\ti (0,\infty)\\[2mm]
\phantom{\prt_{t}u_{m}|u_{m}|^{q-1}}
\lim_{t\to 0}u_{m}(x,t)=\infty\quad\text{locally uniformly on }\Gw^{m}.
\EA\right.\end{equation}
Such a $u_{m}$ is the increasing limit as $k\to\infty$ of the solutions $u_{m,k}$ of the same equation, with same boundary data and initial value equal to $k$. Since $\overline\Gw^{m}\subset\Gw^{m+1}$, $u_{m}<u_{m+1}$. We extend  $u_{m}$
 by zero outside $\Gw^{m}$ and the limit of the sequence $\{u_{m}\}$, when $m\to\infty$ is a positive solution of (\ref{E1}) in 
$Q_{\infty}^{\Gw}$. We denote it by $ u_{\Gw}$. The next result is similar to \rcor{ord}, although the proof is much simpler. 
\bcor{ord+} Assume $\Gw^{1}\subset \Gw^{2}\subset\BBR^N$ are open domains, then $ u_{\Gw^{1}}\leq  u_{\Gw^{2}}$. Furthermore, if  $\Gw=\cup\Gw^{n}$ where $\Gw^{n}\subset\Gw^{n+1}$, then
\begin{equation}\label{comp+}
\lim_{n\to\infty} u_{\Gw^{n}}= u_{\Gw},
\end{equation} 
locally uniformly in $Q_{\infty}^\Gw$.\es


\bprop{equ}There holds $u_{\Gw}=\underline u_{\Gw}$.
\es
\Proof For any $m,k>0$,  $\tilde u_{m,k}$, the extension of  $u_{m,k}$ by zero in $Q_{\infty}^{{\Gw^{m}}^c}$ is a subsolution, thus it is dominated by $\underline u_{\Gw}$. Letting successively $k\to\infty$ and $m\to\infty$ implies $u_{\Gw}\leq\underline u_{\Gw}$.
In order to prove the reverse inequality, we consider an increasing sequence  $\{u_{\ell}\}\subset \CI(Q_{\infty}^\Gw)$ converging to 
$\underline u_{\Gw}$ locally uniformly in $Q_{\infty}^{\Gw}$. If $\Gw$ is bounded there exists a bounded sequence $\{u_{\ell,0,k}\}$ which converges to $u_{\ell}(.,0)=u_{\ell,0}$ in $L^2(\Gw)$ and $S^{\prt J_{\Gw}}(.)[u_{\ell,0,k}]\to S^{\prt J_{\Gw}}(.)[u_{\ell,0}]$ in $L^\infty(0,\infty;L^2(\Gw))$. Therefore 

\begin{equation}\label{equ1}
S^{\prt J_{\Gw}}(.)[u_{\ell,0,k}]\leq u_{\Gw}\Longrightarrow 
S^{\prt J_{\Gw}}(.)[u_{\ell,0}]\leq u_{\Gw}\Longrightarrow\underline u_{\Gw}\leq u_{\Gw}.\end{equation}
Next, if $\Gw$ is unbounded, $\Gw=\cup\Gw^{n}$, with $\Gw^{n}\subset\Gw^{n+1}$ are bounded, we have 
$$\lim_{n\to\infty}\underline u_{\Gw^{n}}=\underline u_{\Gw}$$
and
$$\lim_{n\to\infty} u_{\Gw^{n}}= u_{\Gw}$$
by \rcor{ord} and  \rcor{ord+}. Since $\underline u_{\Gw^{n}}=u_{\Gw^{n}}$ from the first part of the proof, the result follows.\qeda
\medskip

\noindent\Remark By construction $u_{\Gw}$ is dominated by any positive solution of (\ref{Equ})  which satisfies the initial blow-up condition locally uniformly in $\Gw$. Therefore, $u_{\Gw}=\underline u_{\Gw}$ is the {\it minimal solution} with initial blow-up.\medskip

If $\Gw$ has the minimal regularity which allows the Dirichlet problem to be solved by any continuous function $g$ given on $\prt\Gw\ti[0,\infty)$, we can consider another construction of the maximal solution of (\ref{MS0}) in $Q_{\infty}^{\Gw}$. The needed assumption on $\prt\Gw$ is known as the {\it parabolic  Wiener criterion} \cite{Zi} (abr. PWC). 

\bdef{J(Q)}If $\prt\Gw$ is compact and satisfies PWC, we denote by $\CJ_{Q_{\infty}^{\Gw}}$ the set of $v\in C((0,\infty)\ti\overline\Gw)\cap C^{2,1}(Q_{\infty}^{\Gw})$ satisfying  (\ref{MS0}).
\es
\bth{max-bis}Assume $q>1$ and $\Gw$ satisfies PWC. Then  $\CJ_{Q_{\infty}^{\Gw}}$ admits a maximal element 
$\overline u_{\Gw}$.
\es
\Proof {\it Step 1- Construction.} We shall directly assume that $\Gw$ is unbounded, the bounded case being a simple adaptation of our construction. We suppose $\Gw^c\subset B_{R_{0}}$, and for $n>R_{0}$ set $\Gw_{n}=\Gw\cap B_{n}$. The construction of $u_{n}$ is standard: for $k\in\BBN_{*}$ we denote by $v^*_{k}=v^*_{n,k}$ the solution of (\ref{I1}). \rlemma {sup} is valid for $v^*_{k}$. Notice that uniqueness follows from the maximum principle. When $k\to\infty$ the sequence $\{v_{k}\}$
increases and converges to a solution $u_{n}$ of (\ref{Equ}) in $Q_{\Gw_{n}}$. Because the exterior boundary of $\Gw_{n}$ is smooth, the standard equi-continuity of the sequence of solutions applies, thus
$u_{n}(x,t)=0$ for all $(x,t)$ s.t. $|x|=n$ and $t>0$. In order to see that $u_{n}(x,t)=0$ for all $(x,t)$ s.t. $x\in\prt\Gw$ and $t>0$, we see that $u_{n}(x,t)\leq \phi_\gt(x,t)$ on $(\gt,\infty)\ti\Gw_{n}$, where
\begin{equation}\label{bis1}\left\{\BA {l}
\prt_{t}\phi_\gt-\Gd\phi_\gt+\phi^{\gt\,q}=0\quad\text{in }Q_{\infty}^{\Gw}\\
\phantom{----}
\phi_\gt(x,\gt)=\phi_{q}(\gt)\quad\text{in }\Gw\\\phantom{----}
\phi_\gt(x,t)=0\quad\text{in }\prt\Gw\ti [\gt,\infty)
\EA\right.\end{equation}
Such a solution exists because of PWC assumption. Since $v^*_{n,k}$ is an increasing function of $n$ (provided the  solution is extended by $0$ outside $\Gw_{n}$) and $k$, there holds $\tilde u_{n}\leq u_{n+1} $ in $\Gw^{n+1}$. If we set 
$$\overline u_{\Gw}=\lim_{n\to\infty}\tilde u_{n},$$
 then $\overline u_{\Gw}\leq \phi_\gt$ for any $\gt>0$. Clearly $\overline u_{\Gw}$ is a solution of (\ref {Equ}) in $Q_{\infty}^{\Gw}$. This implies that $\overline u_{\Gw}$ is continuous up to $\prt\Gw\ti (0,\infty)$, with zero boundary value. Thus it belongs to $\CJ_{Q_{\infty}^{\Gw}}$. \smallskip
 
 \noindent {\it Step 2- Comparison.} In order to compare $\overline u_{\Gw}$ to any other $u\in\CJ_{Q_{\infty}^{\Gw}}$, for $R>R_{0}$ we set $v_{R,\gt}(x,t)=\overline u_{\Gw}(x,t-\gt)+W_{R}(x)$, where $W_{R}$ is the maximal solution of (\ref{maxell}) in $B_{R}$. The function $(u-v_{R,\gt})_{+}$ is a subsolution of (\ref{Equ}) in $\Gw\cap B_{R}\ti (\gt,\infty)$. It vanishes in a neighborhood on $\prt(\Gw\cap B_{R})\ti (\gt,\infty)$ and of $\Gw\cap B_{R}\ti \{\gt\}$. Thus it is identically zero. If we let $R\to\infty$ in the inequality $u\leq v_{R,\gt}$ and $\gt\to 0$, we derive $u\leq \overline u_{\Gw}$, which is the claim.
\qeda\medskip

\bprop{reg} Under the assumptions of \rth{max-bis}, $ \overline u_{\Gw}\in \CI(Q_{\infty}^{\Gw })$ if $\Gw$ is bounded and $ \overline u_{\Gw}\in \CI(Q_{\infty}^{\Gw_{loc}})$ if $\Gw$ is not bounded.\es
\Proof {\it Case 1: $\Gw$  bounded}.
Let $\Gw^{n}$ be a sequence of smooth domains such that 
$$\Gw^{n}\subset\overline{\Gw^{n}}\subset\Gw^{n+1}\subset\Gw$$
and $\cup_{n}\Gw^{n}=\Gw$. For $\gt>0$, let $u_{n,\gt}$ be the solution of 
\begin{equation}\label{A1}\left\{\BA {l}
\prt_{t}u_{n,\gt}-\Gd u_{n,\gt}+u_{n,\gt}^q=0\quad\text{in }\Gw^{n}\ti (\gt,\infty)\\[2mm]
\phantom{------}
u_{n,\gt}(.,\gt)=\overline u_{\Gw}(.,\gt)\quad\text{in }\Gw^{n}\\[2mm]\phantom{------}
u_{n,\gt}(x,t)=0\quad\text{in }\prt\Gw^{n}\ti [\gt,\infty)
\EA\right.\end{equation}
Because $\overline u_{\Gw}(.,\gt)\in C^2(\overline\Gw^{n})$, $u_{n,\gt}\in C^{2,1}(\overline \Gw^{n}\ti [\gt,\infty))$. By the maximum principle,
\begin{equation}\label{A2}
0\leq \overline u_{\Gw}(.,t)-u_{n,\gt}(.,t)\leq \max\{\overline u_{\Gw}(x,s):(x,s)\in \prt\Gw^{n}\ti [\gt,t]\}
\end{equation}
for any $t>\gt$. Because $\overline u_{\Gw}$ vanishes on $\prt\Gw\ti [\gt,t]$, we derive
\begin{equation}\label{A2'}\lim_{n\to\infty}\tilde u_{n,\gt}= \overline u_{\Gw}
\end{equation}
uniformly on $\overline \Gw\ti [\gt,t]$ for any $t\geq\gt$, where $\tilde u_{n,\gt}$ is the extension of $u_{n,\gt}$ by zero outside $\Gw_{n}$. Applying (\ref{I6}) and (\ref{I10}) with $\eta=1$ to $\tilde u_{n,\gt}$ in $\Gw$ yields to
\begin{equation}\label{A3}\BA {l}2^{-1}\myint{\Gw}{}\tilde u_{n,\gt}^2(.,t)dx+\myint{\gt}{t}\myint{\Gw}{}\left(\abs{\nabla \tilde u_{n,\gt}}^2 +\tilde u_{n,\gt}^{q+1} \right) dxds\leq C(t+1)\gt^{-2/(q-1)}.
\EA\end{equation}
and
\begin{equation}\label{A4}\BA {l}
2^{-1}\myint{\gt}{t}\myint{\Gw}{}(s-\gt)(\prt_{s}\tilde u_{n,\gt})^2 dxds+(t-\gt)\myint{\Gw}{}\left(\myfrac{|\nabla{\tilde u_{n,\gt}}|^2}{2}+
\myfrac{\tilde u_{n,\gt}^{q+1}}{q+1}\right)(t,.)dx
\leq C(t-\gt)t\gt^{-2/(q-1)}.
\EA\end{equation}
Letting $n\to\infty$ and using (\ref{A2'}) yields to
\begin{equation}\label{A5}\BA {l}2^{-1}\myint{\Gw}{}\overline u_{\Gw}^2(.,t)dx+\myint{\gt}{t}\myint{\Gw}{}\left(\abs{\nabla\overline u_{\Gw}}^2 +\overline u_{\Gw}^{q+1} \right) dxds\leq C(t+1)\gt^{-2/(q-1)}.
\EA\end{equation}
and
\begin{equation}\label{A'5}\BA {l}
2^{-1}\myint{\gt}{t}\myint{\Gw}{}(s-\gt)(\prt_{s}\overline u_{\Gw})^2 dxds+(t-\gt)\myint{\Gw}{}\left(\myfrac{|\nabla{\overline u_{\Gw}}|^2}{2}+
\myfrac{\overline u_{\Gw}^{q+1}}{q+1}\right)(t,.)dx
\leq C(t-\gt)t\gt^{-2/(q-1)}.
\EA\end{equation}
Since $L^2(\gt,t;H^1_{0}(\Gw))$ is a closed subspace of $L^2(\gt,t;H^1(\Gw))$, for any $0<\gt <t$, $\overline u_{\Gw}\in L^2_{loc}(0,\infty;H^1_{0}(\Gw))$. Furthermore $\prt_{s}\overline u_{\Gw}\in L^2_{loc}(0,\infty;L^2(\Gw))$. Because $\overline u_{\Gw}$ satisfies (\ref{Equ}), it implies $ \overline u_{\Gw}\in \CI(Q_{\infty}^{\Gw})$. \smallskip

\noindent {\it Case 2: $\Gw$  unbounded}. We assume that $\Gw^c\subset B_{R_{0}}$. We consider a sequence of smooth unbounded domains $\{\Gw^n\}\subset\Gw$ ($n>1$) such that 
$\sup\{\dist(x,\Gw^c):x\in\prt\Gw^{n}\}<1/n$ as $n\to\infty$, thus
$\cup_{n}\Gw^n=\Gw$. For $m>R_{0}$ we set $\Gw_{m}^n=\Gw^n\cap B_{m}$. Therefore
$\Gw_{m}^n\subset\overline{\Gw_{m}^n}\subset\Gw_{m+1}^{n+1}$ and $\cup_{n,m}\Gw_{m}^n=\Gw$. For $\gt>0$, let $u=u_{m,n,\gt}$ be the solution of
\begin{equation}\label{A6}\left\{\BA {l}
\prt_{t}u-\Gd u+u^q=0\quad\text{in }\Gw_{m}^n\ti (\gt,\infty)\\[2mm]
\phantom{----}
u(.,\gt)=\overline u_{\Gw}(.,\gt)\quad\text{in }\Gw_{m}^n\\[2mm]\phantom{----}
u(x,t)=0\quad\text{in }\prt\Gw^n\ti [\gt,\infty)
\\[2mm]\phantom{----}
u(.,\gt)=\overline u_{\Gw}(.,\gt)\quad\text{in }\prt B_{m}\ti (\gt,\infty).\\[2mm]
\EA\right.\end{equation}
By the maximum principle, 
\begin{equation}\label{A7}
0\leq \overline u_{\Gw}(.,t)-u_{m,n,\gt}(.,t)\leq \max\{\overline u_{\Gw}(x,s):(x,s)\in \prt\Gw^{n}\ti [\gt,t]\}\to 0,
\end{equation}
as $n\to 0$. Next we extend $u_{m,n,\gt}$ by zero in $\Gw\setminus\Gw_{n}$ and apply (\ref{I6})-(\ref{I10}) with $\eta$ as in \rth{max} and $m>2R$. We get, with 
$\Gw_R=\Gw\cap B_{R}$,
\begin{equation}\label{A8}\BA {l}
2^{-1}\myint{\Gw_R}{} u_{m,n,\gt}^2(.,t)dx+\myint{\gt}{t}\myint{\Gw_R}{}\left(2^{-1}\abs{\nabla  u_{m,n,\gt}}^2 +u_{m,n,\gt}^{q+1} \right) dxds\leq CR^N(t+1)\gt^{-2/(q-1)},
\EA\end{equation}
and
\begin{equation}\label{A9}\BA {l}
2^{-1}\myint{\gt}{t}\myint{\Gw_R}{}(s-\gt)
\abs{\myfrac{du_{m,n,\gt}}{dt}}^2 dxds+(t-\gt)\myint{\Gw_R}{}
\left(\myfrac{|\nabla{u_{m,n,\gt}|^2}}{2}+\myfrac{u_{m,n,\gt}^{q+1}}{q+1}\right)dx\\[4mm]
\phantom{----------------------}
\leq 8C(2R)^N(t-\gt)t\gt^{-2/(q-1)}
\EA\end{equation}
We let successively $m\to\infty$ and $n\to\infty$ and derive by Fatou's lemma and (\ref{A7}) that inequalities (\ref{A8}) and (\ref{A9}) still hold with $\overline u_{\Gw}$ instead of $u_{m,n,\gt}$. If we denote by $\tilde H_{0}^1(\Gw_{R})$ the closure of the space of $C^\infty(\overline{\Gw_{R}})$ functions which vanish in a neighborhood on $\prt\Gw$, then (\ref{A8}) is an estimate in $L^2(\gt,t;\tilde H_{0}^1(\Gw_{R}))$ which is a closed subspace of $L^2(\gt,t; H^1(\Gw_{R}))$. Therefore 
$\overline u_{\Gw}\in L_{loc}^2(0,\infty;\tilde H_{0}^1(\Gw_{R}))$. Using (\ref{A9}) and equation (\ref{Equ}) we conclude that $ \overline u_{\Gw}\in \CI(Q_{\infty}^{\Gw_{loc}})$.\qeda
\medskip

We end this section with a comparison result between 
$\underline u_{\Gw}$ and $\overline u_{\Gw}$.

\bth{equality} Assume $q>1$ and $\Gw$ satisfies PWC. Then $\underline u_{\Gw}=\overline u_{\Gw}$.
\es
\Proof By \rprop{equ} and \rth{max-bis}-Step 2, $\underline u_{\Gw}\leq \overline u_{\Gw}$. If $\Gw$ is bounded, we can compare $\underline u_{\Gw}(.,.)$ and 
$ \overline u_{\Gw}(.+\gt,.)$ on $\Gw\ti (0,\infty)$. Since 
$\underline u_{\Gw}$, the least upper bound of $\CI(Q^\Gw_{\infty})$ belongs to 
$\CI(Q^\Gw_{\infty})$, and $ \overline u_{\Gw}(.+\gt,.)\in \CI(Q^\Gw_{\infty})$ we derive $ \overline u_{\Gw}(.+\gt,.)\leq \underline u_{\Gw}(.,.)$, from which follows
$\overline u_{\Gw}\leq \underline u_{\Gw}$. Next, if $\Gw$ is not bounded, we can proceed as in the proof of \rth{max} by comparing $\underline u_{\Gw}(.,.)+W_{R}$ and  $ \overline u_{\Gw}(.+\gt,.)$ on $\Gw_{R}\ti (0,\infty)$, where $W_{R}$ is defined in                  (\ref{maxell}). Because $\left(\overline u_{\Gw}(.+\gt,.)-\underline u_{\Gw}(.,.)-W_{R}(.)\right)_{+}$ is a subsolution of (\ref{Equ}) in $Q^{\Gw_{R}}_{\infty}$ which vanishes at 
$t=0$ and near $\prt\Gw_{R}\ti (0,\infty)$; it follows $\overline u_{\Gw}(.+\gt,.)\leq \underline u_{\Gw}(.,.)+W_{R}(.)$. Letting $R\to\infty$ and $\gt\to 0$ completes the proof.\qeda

\section{Uniqueness of large solutions}

\bdef {LS}Let $q>1$ and $\Gw\subset\BBR^N$ be any domain. A positive function $u\in C^{2,1}(Q_{\infty}^{\Gw})$ of (\ref{Equ}) is a large initial solution if it satisfies
\begin{equation}\label{L1}
\lim_{t\to 0}u(x,t)=\infty\quad\forall x\in\Gw,
\end{equation} 
uniformly on any compact subset of $\Gw$.
\es
We start with the following lemma

\blemma{below1} Assume $u\in C^{2,1}(Q_{\infty}^{\Gw})$ is a large solution of (\ref{Equ}), then for any open subset $G$ such that 
$\overline G\subset\Gw$, there holds
\begin{equation}\label{L2}
\lim_{t\to 0}t^{1/(q-1)}u(x,t)=c_{q}:=\left(\myfrac{1}{q-1}\right)^{1/(q-1)}
\quad\text {uniformly in }G.
\end{equation}
\es
\Proof By compactness, it is sufficient to prove the result when $G=B_{\gr}$ and $\overline B_{\gr}\subset B_{\gr'}\subset\Gw.$ Let $\gt>0$; by comparison,  
$u(x,t)\geq u_{B_{\gr'}}(x,t+\gt)$ for any $(x,t)\in Q_{\infty}^{\Gw}$. Letting $\gt\to 0$ yields to $u\geq u_{B_{\gr'}}$. Next for $\gt>0$,
$$\phi_{q}(t+\gt)\leq u_{B_{\gr'}}(x,t)+u_{B^c_{\gr'}}(x,t)+W_{R}(x)\quad\forall (x,t)\in Q_{\infty}^{\BBR^N}.
$$
Similarly
$$\max\{u_{B_{\gr'}}(x,t+\gt),u_{B^c_{\gr'}}(x,t+\gt)\}\leq \phi_{q}(t)+W_{R}(x)
\quad\forall (x,t)\in Q_{\infty}^{\BBR^N}.
$$
Letting $R\to \infty$ and $\gt\to 0$,
$$\max\{u_{B_{\gr'}},u_{B^c_{\gr'}}\}\leq \phi_{q}\leq u_{B_{\gr'}}+u_{B^c_{\gr'}}\quad\text {in }Q_{\infty}^{\BBR^N}.$$
 For symmetry reasons, $x\mapsto u_{B^c_{\gr'}}(x,t)$ is radially increasing for any $t>0$, thus, for any $\gr<\gr'$ and $T>0$, there exists $C_{\gr,T}>0$ such that 
$$u_{B^c_{\gr'}}(x,t)\leq C_{\gr,T}\quad\forall (x,t)\in B_{\gr}\ti [0,T].$$
Therefore
$$\lim_{t\to 0}t^{1/(q-1)}u_{B_{\gr'}}(x,t)=c_{q}
\quad\text{uniformly on }B_{\gr}.
$$
Because 
$$u_{B_{\gr'}}(x,t)\leq u(x,t)\leq \phi_{q}(t)\quad\forall (x,t)\in Q_{\infty}^{\Gw},
$$
(\ref{L2}) follows.\qeda\medskip

As an immediate consequence of \rlemma {below1} and (\ref{comp}), we obtain

\bprop{M=LS} Assume $q>1$ and $\prt\Gw$ is compact. Then $u_{\Gw}$ is a large solution.
\es

We start with the following uniqueness result

\bprop{starsha} Assume $q>1$, $\Gw$ satisfies PWC, $\prt\Gw$ is bounded, and either $\Gw$ or $\Gw^c$ is strictly starshaped with respect to some point. Then $\overline u_{\Gw}$ is the unique large solution belonging to $\CJ(Q_{\infty}^{\Gw})$.
\es
\Proof Without loss of generality, we can suppose that either $\Gw$ or $\Gw^c$ is strictly starshaped with respect to $0$. By \rth{max-bis}, $\overline u_{\Gw}$ exists and, by (\ref{comp}) and \rlemma {below1},  it is a large solution. Let $u\in \CJ(Q_{\infty}^{\Gw})$ be another large solution. Clearly $u\leq \overline u_{\Gw}$. If $\Gw$ is starshaped, then for $k>1$, the function $u_{k}(x,t):=k^{2/(q-1)}u(kx,k^2t)$ is a solution in $Q_{\Gw_{k}}$, with $\Gw_{k}:=k^{-1}\Gw$. Clearly it is a large solution and it belongs to $\CJ(Q_{\Gw_{k}})$. For $\gt\in (0,1)$, set $u_{k,\gt}(x,t)=u_{k}(x,t-\gt)$. Because $\prt\Gw$ is compact, 
$$\lim_{k\downarrow 1}d_{H}(\prt\Gw,\prt\Gw_{k})=0,$$
where $d_{H}$ denotes the Hausdorff distance between compact sets. By assumption 
$\overline u_{\Gw}\in C([\gt,\infty)\ti\overline\Gw)$ vanishes on $[\gt,\infty)\ti\prt\Gw$, thus, for any $\ge>0$, there exists $k_{0}>1$ such that for any 
$$k\in (1,k_{0}]\Longrightarrow \sup\{\overline u_{\Gw}(x,t):(x,t)\in [\gt,1]\ti \prt \Gw_{k}\}\leq \ge.$$
Since $u_{k,\gt}+\ge$ is a super solution in $Q_{\Gw_{k}}$ which dominates $\overline u_{\Gw}$ on $[\gt,1]\ti \prt \Gw_{k}$ and at $t=\gt$, it follows that 
$u_{k,\gt}+\ge\geq \overline u_{\Gw}$ in $(\gt,1]\ti \Gw_{k}$. Letting successively 
$k\to 1$, $\gt\to 0$ and using the fact that $\ge$ is arbitrary, yields to $u\geq \overline u_{\Gw}$ in $(0,1]\ti \Gw$ and thus in $Q_{\infty}^{\Gw}$. If $\Gw^c$ is starshaped, then the same construction holds provided we take $k<1$ and use the fact that, for $R>0$ large enough, $u_{k,\gt}+\ge+W_{R}$ is a super solution in $Q_{\Gw_{k}\cap B_{R}}$ which dominates $\overline u_{\Gw}$ on $[\gt,1]\ti \prt \Gw_{k}\cap B_{R}$ and at $t=\gt$. Letting successively $R\to\infty$, 
$k\to 1$, $\gt\to 0$ and $\ge\to 0$ yields to $u\geq \overline u_{\Gw}$\qeda\medskip

As a consequence of Section 2, we have the more complete uniqueness theorem

\bth{Uniq1} Assume $q>1$, $\Gw\subset\BBR^N$ is a domain with a bounded boundary $\prt\Gw$ satisfying PWC. Then for any $f\in C(\prt\Gw\ti[ 0,\infty))$, $f\geq 0$, there exists a unique positive function $u=\overline u_{\Gw,f}\in C(\overline \Gw\ti(0,\infty))\cap C^{2,1}(Q_{\infty}^\Gw)$ satisfying
\begin{equation}\label{CD}\left\{\BA {l}
\prt_{t}u-\Gd u+|u|^{q-1}u=0\quad\text{in }Q_{\infty}^{\Gw}\\[2mm]
\phantom{\prt_{t}u-\Gd u+|u|^{q-1}}
u=f\quad\text{in }\prt\Gw\ti (0,\infty)\\[2mm]
\phantom{,u^{q-1}}
\lim_{t\to 0}u(x,t)=\infty\quad\text{locally uniformly on }\Gw.
\EA\right.\end{equation}
\es
\Proof {\it Step 1: Existence.} It is a simple adaptation of the proof of \rth{max-bis}. For $k,\gt>0$, we denote by $u=u_{k,\gt,f}$ the solution of 
\begin{equation}\label{CD2}\left\{\BA {l}
\prt_{t}u-\Gd u+|u|^{q-1}u=0\quad\text{in }\Gw\ti (\gt,\infty)\\[2mm]
\phantom{\prt_{t}u-\Gd u+|u|^{q-1}}
u=f\quad\text{in }\prt\Gw\ti (\gt,\infty)\\[2mm]
\phantom{,u^{q-1}---}
u(x,\gt)=k\quad\text{on }\Gw.
\EA\right.\end{equation}
Notice that $u_{k,\gt,f}$ is bounded from above by $\overline u_{\Gw}(.,.-\gt)+v_{f,\gt}$, where $v_{f,\gt}=v$ solves
\begin{equation}\label{CD3}\left\{\BA {l}
\prt_{t}v-\Gd v+|v|^{q-1}v=0\quad\text{in }\Gw\ti (\gt,\infty)\\[2mm]
\phantom{\prt_{t}v-\Gd v+|v|^{q-1}}
v=f\quad\text{in }\prt\Gw\ti (\gt,\infty)\\[2mm]
\phantom{,u^{q-1}---}
v(x,\gt)=0\quad\text{on }\Gw.
\EA\right.\end{equation}
If we let $k\to\infty$ we obtain a solution $u_{\infty,\gt,f}$ of the same problem except that the condition at $t=\gt$ becomes $\lim_{t\to\gt}u(x,t)=\infty$, locally uniformly for $x\in\Gw$. Clearly $u_{\infty,\gt,f}$ dominates in $\Gw\ti (\gt,\infty)$ the restriction to this set of any $u \in C(\overline \Gw\ti\infty))\cap C^{2,1}(Q_{\infty}^\Gw)$ solution of (\ref{CD}), in particular $\overline u_{\Gw}$. Therefore $u_{\infty,\gt,f}\geq u_{\infty,\gt',f}$ in 
$\Gw\ti (\gt,\infty)$ for any $0<\gt'<\gt$. When $\gt\to 0$, $u_{\infty,\gt,f}$ converges to a function $\overline u_{f}$ which satisfies the lateral boundary condition $\overline u_{\Gw,f}=f$. Therefore $\overline u_{\Gw,f}$ satisfies (\ref{CD}).\smallskip

\noindent {\it Step 2: Uniqueness.} Assume that there exists another positive function $u:=u_{f}\in C(\overline \Gw\ti(0,\infty))\cap C^{2,1}(Q_{\infty}^\Gw)$ solution of (\ref{CD}). Then $u_{f}<\overline u_{\Gw,f}$. For $\gt>0$, consider the solution $v:=v_{\gt}$ of 
\begin{equation}\label{CD4}\left\{\BA {l}
\prt_{t}v-\Gd v+|v|^{q-1}v=0\quad\text{in }\Gw\ti (\gt,\infty)\\[2mm]
\phantom{\prt_{t}v-\Gd v+|v|^{q-1}}
v=0\quad\text{in }\prt\Gw\ti (\gt,\infty)\\[2mm]
\phantom{,u^{q-1}---}
v(x,\gt)=u_{f}(x,\gt)\quad\text{on }\Gw.
\EA\right.\end{equation}
Then $v_{\gt}\leq u_{f}$ in $\Gw\ti (\gt,\infty)$. In the same way, we construct a solution $v:\tilde v_{\gt}$ of the same problem (\ref{CD4}) except that the condition at $t=\gt$ is now $v(x,\gt)=\overline u_{\Gw,f}(x,\gt)$ for all $x\in\Gw$.
Furthermore $v_{\gt}\leq \tilde v_{\gt}\leq \overline u_{\Gw,f}$. Next we adapt a method introduced in \cite {MV1}, \cite {MV2} in a different context. We denote
\begin{equation}\label{CD5}
Z_{f}=\overline u_{\Gw,f}-u_{f}\quad\text{and }\;Z_{0,\gt}=\tilde v_{\gt}-v_{\gt},
\end{equation}
and, for $(r,s)\in \BBR_{+}^2$,
$$h(r,s)=\left\{\BA {l}\myfrac{r^q-s^q}{r-s}\quad\text{if  }r\neq s\\
0\qquad\qquad\text{if }r= s.
\EA\right.$$
Since $r\mapsto r^q$ is convex on $\BBR_{+}$, there holds
$$
\left\{\BA {l}r_{0}\geq s_{0},\,r_{1}\geq s_{1}\\
r_{1}\geq r_{0},\,s_{1}\geq s_{0}
\EA\right.\Longrightarrow h(r_{1},s_{1})\geq h(r_{0},s_{0}).
$$
This implies
\begin{equation}\label{CD6}h(u_{\Gw,f},u_{f})\geq h(\tilde v_{\gt},v_{\gt})\quad\text{in }
\Gw\ti [\gt,\infty).
\end{equation}
Next we write
\begin{equation}\label{CD7}\BA {l}
0=\prt_{t}(Z_{f}-Z_{0,\gt})-\Gd(Z_{f}-Z_{0,\gt})+\overline u^q_{\Gw,f}-u^q_{f}-(\tilde v^q_{\gt}-v^q_{\gt})\\[2mm]\phantom{0}=
\prt_{t}(Z_{f}-Z_{0,\gt})-\Gd(Z_{f}-Z_{0,\gt})+h(\overline u_{\Gw,f},u_{f})Z_{f}-h(\tilde v_{\gt},v_{\gt})Z_{0,\gt}.
\EA\end{equation}
Combining (\ref{CD6}), (\ref{CD7}) with the positivity of  $Z_{f}$ and $Z_{0,\gt}$, we derive
\begin{equation}\label{CD8}
\prt_{t}(Z_{f}-Z_{0,\gt})-\Gd(Z_{f}-Z_{0,\gt})+h(\overline u_{\Gw,f},u_{f})(Z_{f}-Z_{0,\gt})\leq 0,
\end{equation}
in $\Gw\ti (\gt,\infty)$. On $\prt\Gw\ti [\gt,\infty)$ there holds
$Z_{f}-Z_{0,\gt}=f-f=0$. Furthermore, at at $t=\gt$,
$Z_{f}(x,\gt)-Z_{0,\gt}(x,\gt)=\overline u_{\Gw,f}(x,t)-u_{f}(x,\gt)-\overline u_{\Gw,f}(x,t)+u_{f}(x,\gt)=0$. By the maximum principle, it follows $Z_{f}\leq Z_{0,\gt}$ in $\Gw\ti [\gt,\infty)$. 
Since $\gt>\gt'>0$ implies $v_{\gt}(x,\gt)=u_{f}(x,\gt)\geq v_{\gt'}(x,\gt)$ and $\tilde v_{\gt}(x,\gt)=\overline u_{\Gw,f}(x,\gt)\geq \tilde v_{\gt'}(x,\gt)$, the sequences $\{v_{\gt}\}$ and $\tilde v_{\gt}$ converge to some functions $\{v_{0}\}$ and $\tilde v_{0}$ which belong to $C(\overline \Gw\ti(0,\infty))\cap C^{2,1}(Q_{\infty}^\Gw)$ and satisfy (\ref{CD}) with $f=0$ on $\prt\Gw\ti (0,\infty)$. Furthermore
\begin{equation}\label{CD9}
\overline u_{\Gw,f}-u_{f}\leq \tilde v_{0}-v_{0}.
\end{equation}
Since 
$\overline u_{\Gw,f}\geq \overline u_{\Gw}$, $\tilde v_{0}\geq \overline u_{\Gw}$, which implies that $\tilde v_{0}=\overline u_{\Gw}$ by the maximality of $\overline u_{\Gw}$. If $\Gw'$ is any smooth bounded open subset such that $\overline\Gw'\subset\Gw$ there holds by an easy approximation argument $v_{0}\geq u_{\Gw'}$ in $\Gw'\ti (0,\infty)$. Therefore $v_{0}\geq u_{\Gw}=\underline  u_{\Gw}=\overline  u_{\Gw}$, by \rprop{equ} and  \rth{equality}. Applying again \rth{equality} we derive that the right-hand side of (\ref{CD9}) is zero, which yields to $\overline u_{\Gw,f}=u_{f}$\qeda

 \end{document}